\documentclass[12pt]{report}
\usepackage{amssymb}
\usepackage{amsmath}

\begin{document}

\vspace{.2in}\parindent=0mm

\begin{center}

{\bf\Large {Construction of Biorthogonal Wavelet Packets on\\\parindent=0mm \vspace{.1in} Local Fields  of Positive Characteristic }}

\vspace{.6in}\parindent=0mm

  {\bf{ F. A. Shah$^{*}$ and M. Y. Bhat$^{**}$  }}

\parindent=0mm \vspace{.2in}
{{\it $^{*}$Department of  Mathematics,  University of Kashmir\\
 South Campus, Anantnag-192101\\
 Jammu and Kashmir, India\\
 $\text{fashah79@gmail.com}$}}

\parindent=0mm \vspace{.1in}
{{\it $^{**}$Department of  Mathematics,  Central University of Jammu\\
Jammu-180011, Jammu and Kashmir, India\\
 $\text{gyounusg@gmail.com}$}}

\end {center}

\parindent=0mm \vspace{.1in}

 {\bf{Abstract:}} Orthogonal wavelet packets lack symmetry which is a much desired property in image and signal processing. The biorthogonal wavelet packets achieve symmetry where the orthogonality is replaced by the biorthogonality. In the present paper, we construct biorthogonal wavelet packets on local fields of positive characteristic and investigate their properties by means of the Fourier transforms. We also show how to obtain several new Riesz bases of the space $L^2(K)$ by constructing a series of subspaces of these wavelet packets. Finally, we provide the algorithms for the decomposition and reconstruction using these biorthogonal wavelet packets.

\parindent=0mm \vspace{.2in}


\parindent=0mm \vspace{.2in}


\parindent=0mm \vspace{.2in}

{\bf{1. Introduction}}

\parindent=8mm\vspace{.2in}
A field $ K$ equipped with a topology is called a {\it local field} if both the additive $K^+$ and multiplicative groups $K^*$ of $ K$ are locally compact Abelian groups. The local fields are essentially of two types: zero and positive characteristic (excluding the connected local fields $\mathbb R$ and $\mathbb C$). Examples of local fields of characteristic zero include the $p$-adic field $\mathbb Q_p$ where as local fields of positive characteristic are the Cantor dyadic group and the Vilenkin $p$-groups. Even though the structures and metrics of local fields of zero and positive characteristics are similar, their wavelet and multiresolution analysis theory are quite different. In recent  years, local fields have attracted the attention of several mathematicians, and have found innumerable applications not only to number theory but also to representation theory, division algebras, quadratic forms and algebraic geometry. As a result, local fields are now consolidated as part of the standard repertoire of contemporary mathematics. For more about local fields and their applications, we refer to the monographs [15, 21].

\parindent=8mm \vspace{.2in}
In recent years, there has been a considerable interest in the problem of constructing wavelet bases on various spaces other than $\mathbb R$, such as abstract Hilbert spaces [20], locally compact Abelian groups [8], Cantor dyadic groups [11], $p$-adic fields [10] and zero-dimensional groups [14]. Recently, R.L. Benedetto and J.J. Benedetto [3] developed a wavelet theory for local fields and related groups. They did not develop the multiresolution analysis (MRA) approach, their method is based on the theory of wavelet sets. The concept of multiresolution analysis on a local field $K$ of positive characteristic was introduced by Jiang {\it et al.} [9]. They pointed out a method for constructing orthogonal wavelets on local field $ K$ with a constant generating sequence and derived necessary and sufficient conditions for a solution of the refinement equation to generate a multiresolution analysis of $L^2( K)$. Subsequently, the tight wavelet frames on local fields were constructed by Li and Jiang  in [13]. They have established  necessary condition and sufficient conditions for tight wavelet frame on local fields  in frequency domain. Shah and Debnath [18] have constructed tight wavelet frames on local field $ K$ of positive characteristic via extension principles. They also provide a sufficient condition for finite number of functions to form a tight wavelet frame and established general principle for constructing tight wavelet frames on local fields. On the other hand, Behera and Jahan [1] have constructed biorthogonal wavelets on a local field $ K$ of positive characteristic and showed that if $\varphi$ and $\tilde \varphi$ are the dual scaling functions assocaited with dual MRA's on local field $K$ of positive characteristic such that their translates are biorthogonal, then the corresponding wavelet families are also biorthogonal.

\parindent=8mm \vspace{.2in}
It is well-known that the classical orthogonal wavelet bases have poor frequency localization. To overcome this disadvantage, Coifman {\it et al.}[7] constructed univariate orthogonal wavelet packets. Well known Daubechies orthogonal wavelets are a special case of wavelet packets. Later on, Chui and Li [5] generalized the concept of orthogonal wavelet packets to the case of non-orthogonal wavelet packets so that they can be applied to the spline wavelets and so on. The introduction of biorthogonal wavelet packets attributes to Cohen and Daubechies [6]. They have also shown that all the wavelet packets, constructed in this way, are not led to Riesz bases for $L^2(\mathbb R)$. Shen [19] generalized the notion of univariate orthogonal wavelet packets to the case of multivariate wavelet packets. Other notable generalizations are the orthogonal version of vector-valued wavelet packets [4],  multiwavelet packets [12], orthogonal and biorthogonal wavelet packets related  to the Walsh polynomials on a positive half-line $\mathbb R^+$ [16,17].

\parindent=8mm \vspace{.2in}
Recently, Behera and Jahan [2] have constructed  orthogonal wavelet packets and wavelet frame packets on local field $ K$ of positive characteristic and show how to construct an orthonormal basis from a Riesz basis. Orthogonal wavelet packets have many desired properties such as compact support, good frequency localization and vanishing moments. However, there is no continuous symmetry which is a much desired property in imaging the compression  and signal processing. To achieve symmetry, several generalizations of scalar orthogonal wavelet packets have been investigated in literature. The biorthogonal wavelet packets achieve symmetry where the orthogonality is replaced by the biorthogonality. Therefore, the objective of this paper is to construct biorthogonal wavelet packets on local fields of positive characteristic and investigate their properties by means of the Fourier transforms and construct several new Riesz bases of space $L^2(K)$. Finally, we establish some algorithms for decomposition and reconstruction using these biorthogonal wavelet packets.

\parindent=8mm \vspace{.2in}
This paper is organized as follows. In Section 2, we discuss some preliminary facts about local fields of positive characteristic and also some  results which are required in the subsequent sections including the definition of an MRA on local fields. In Section 3, we examined some of the properties of the biorthogonal wavelet packets via Fourier transforms. In Section 4, we generate Riesz bases of $L^2(K)$ from these wavelet packets. Section 5, deals with the decomposition and reconstruction algorithms corresponding to these wavelet packets.

\parindent=0mm \vspace{.2in}

\pagestyle{myheadings}

{\bf{2.  Preliminaries on Local Fields }}

\parindent=8mm \vspace{.2in}
Let $K$ be a field and a topological space. Then $K$ is called a {\it local field} if both  $K^+$ and $K^*$ are locally compact Abelian groups, where  $K^+$ and $K^*$ denote the additive and multiplicative groups of $K$, respectively. If $K$ is any field and is endowed with the discrete topology, then $K$ is a local field. Further, if $K$ is connected, then $K$ is either $\mathbb R$ or $\mathbb C$. If $K$ is not connected, then it is totally disconnected. Hence by a local field, we mean a field $K$ which is locally compact, non-discrete and totally disconnected. The $p$-adic fields are examples of local fields. More details are referred to [15, 21]. In the rest of this paper, we use $\mathbb N, \mathbb N_0$ and $\mathbb Z$ to denote the sets of natural, non-negative integers and integers, respectively.

\parindent=8mm \vspace{.2in}

Let $K$ be a fixed local field. Then there is an integer $q = {\mathfrak p}^r$, where $\mathfrak p$ is a fixed prime element of
$K$ and $r$ is a positive integer, and a norm $| . |$ on $K$ such that for all $x \in K$ we have $|x|\ge 0$ and for each $x\in K\setminus \left\{ 0\right\}$ we get $|x| = q^k$ for some integer $k$. This norm is non-Archimedean, that is
$|x+y|\le \max \left\{ |x|, |y|\right\}$ for all $x, y\in K$ and $|x+y|= \max \left\{ |x|, |y|\right\}$ whenever $|x|\not= |y|$.
Let $dx$ be the Haar measure on the locally compact, topological group $(K,+)$. This measure is normalized so that $\int_{\mathfrak D}
dx = 1$, where ${\mathfrak D}= \left\{x \in K: |x| \le 1\right\}$ is the {\it ring of integers} in ${K}$. Define ${\mathfrak B}= \left\{x \in K: |x| < 1\right\}$. The set ${\mathfrak B}$ is called the {\it prime ideal} in $K$. The prime ideal in $K$ is the unique maximal ideal in ${\mathfrak D}$ and hence as result ${\mathfrak B}$ is both principal and prime. Therefore, for such an ideal ${\mathfrak B}$ in ${\mathfrak D}$, we have ${\mathfrak B}= \langle \mathfrak p \rangle=\mathfrak p {\mathfrak D}.$

\parindent=8mm \vspace{.2in}
Let ${\mathfrak D}^*= {\mathfrak D}\setminus {\mathfrak B }=\left\{x\in K: |x|=1   \right\}$. Then, it is easy to verify that ${\mathfrak D}^*$ is a group of units in $K^*$ and if $x\not=0$, then we may write $x=\mathfrak p^k x^\prime, x^\prime\in {\mathfrak D}^*.$ Moreover, each ${\mathfrak B}^k= \mathfrak p^k {\mathfrak D}=\left\{x \in K: |x| < q^{-k}\right\}$ is a compact subgroup of $K^+$ and usually  known as the {\it fractional ideals} of $K^+$ (see [15]). Let ${\cal U}= \left\{a_i \right\}_{ i=0}^{q-1}$ be any fixed full set of coset representatives of ${\mathfrak B}$ in ${\mathfrak D}$, then every element $x\in K$ can be expressed uniquely  as $x=\sum_{\ell=k}^{\infty} c_\ell \mathfrak p^\ell $ with $c_\ell \in {\cal U}.$ Let $\chi$ be a fixed character on $K^+$ that is trivial on ${\mathfrak D}$ but is non-trivial on  ${\mathfrak B}^{-1}$. Therefore, $\chi$ is constant on cosets of ${\mathfrak D}$ so if $y \in {\mathfrak B}^k$, then $\chi_y(x)=\chi(yx), x\in K.$ Suppose that $\chi_u$ is any character on $K^+$, then clearly the restriction $\chi_u|{\mathfrak D}$ is also a character on ${\mathfrak D}$. Therefore, if $\left\{u(n): n\in\mathbb N_0\right\}$ is a complete list of distinct coset representative of ${\mathfrak D}$ in $K^+$, then it is proved in [21] that the set  $\left\{\chi_{u(n)}: n\in\mathbb N_0\right\}$   of distinct characters on ${\mathfrak D}$ is a complete orthonormal system on ${\mathfrak D}$.

\parindent=0mm \vspace{.2in}
The Fourier transform $\hat f$ of a function $f \in L^1(K)\cap L^2(K)$ is defined by\\
$$\hat f(\xi)= \displaystyle \int_K f(x)\overline{ \chi_\xi(x)}dx.\eqno(2.1)$$
It is noted that\\
$$\hat f(\xi)= \displaystyle \int_K f(x) \overline{ \chi_\xi(x)}dx= \displaystyle \int_K f(x)\chi(-\xi x)dx.$$

\parindent=8mm \vspace{.2in}
Furthermore, the properties of Fourier transform on local field are much similar to those of on the real line. In particular Fourier transform is unitary on $L^2(K)$.

\parindent=8mm \vspace{.2in}
Let us now impose a natural order on the sequence $\{u(n)\}_{n=0}^\infty$. Since ${\mathfrak D}/ \mathfrak B \cong GF(q) $ where $ GF(q) $ is a $c$-dimensional vector space over the field $ GF(q) $, we choose a set $\{1=\zeta_0,\zeta_1,\zeta_2,\dots,\zeta_{c-1}\}\subset {\mathfrak D^*}$ such that span$\{\zeta_j\}_{j=0}^{c-1}\cong GF(q)$. Let $\mathbb N_0=\mathbb N\cup \{0\}$. For $n \in \mathbb N_0$ such that $0\leq n<q$, we have
$$n=a_0+a_1{\mathfrak p}+\dots+a_{c-1}{\mathfrak p}^{c-1},\;0\leq a_k<{\mathfrak p},\,k=0,1,\dots,c-1.$$
Define
$$u(n)=(a_0+a_1\zeta_1+\dots+a_{c-1}\zeta_{c-1}){\mathfrak p}^{-1}.\eqno(2.2)$$

\parindent=0mm \vspace{.1in}
For $n \in \mathbb N_{0}$ and $0\leq b_k<q,k=0,1,2,\dots,s$, we write
$$n=b_0+b_1q+b_2q^2+\dots+b_sq^s,$$
such that
$$u(n)=u(b_0)+u(b_1){\mathfrak p}^{-1}+\dots+u(b_s){\mathfrak p}^{-s}.\eqno(2.3)$$

\parindent=8mm \vspace{.1in}
If $r,k\in\mathbb N_{0}\; \text{and}\;0\le s<q^k$, then it follows that
$$u(rq^k+s)=u(r){\mathfrak p}^{-k}+u(s).$$

\parindent=0mm \vspace{.0in}
Further, it is easy to verify that $u(n)=0\;\text{if and only if} \; n=0$ and $\{u(\ell)+u(k):k \in \mathbb N_0\}=\{u(k):k \in \mathbb N_0\}\;\text{for a fixed}\; \ell \in \mathbb N_0.$ From now we will write $\chi_n$ for $\chi_{u(n)}$.

\parindent=8mm \vspace{.2in}
Let the local field $K$ be of characteristic $p>0$ and $\zeta_0,\zeta_1,\zeta_2,\dots,\zeta_{c-1}$ be as above. We define a character $\chi$ on $K$ as follows:
$$\chi(\zeta_\mu)= \left\{
\begin{array}{lcl}
\exp(2\pi i/p),&&\mu=0\;\text{and}\;j=1,\\
1,&&\mu=1,\dots,c-1\;\text\;j \neq 1.
\end{array}
\right. \eqno(2.4)
$$

\parindent=0mm \vspace{.2in}
{\bf{Definition 2.1.}} Let $\{f_n\}_{n\in \mathbb N_0}$ be a sequence of a Hilbert space $H$. Then, $\{f_n\}_{n\in \mathbb N_0}$ is said to form a {\it{Riesz basis}} for $H$ if

\parindent=0mm \vspace{.1in}
(a) $\{f_n\}_{n\in \mathbb N_0}$ is linearly independent, and

\parindent=0mm \vspace{.1in}
(b) there exists constants $A$ and $B$ with $0<A\le B<\infty$ such that
$$A\big|\big|f\big|\big|_2^2\le\displaystyle\sum_{k\in \mathbb N_0}\left|\left<f, f_n\right>\right|^2 \le B\big|\big|f\big|\big|_2^2 \quad \text{for every}\;f \in H.\eqno(2.5)$$

\parindent=8mm \vspace{.0in}
Let us recall the definition of an MRA on local fields of positive characteristics ([9]).

\parindent=0mm \vspace{.2in}
{\bf{Definition 2.2.}} Let $K$ be a local field of positive characteristic $p> 0$ and $ \mathfrak p$ be a prime element of $K$. A multiresolution analysis(MRA) of $L^2(K)$ is a sequence of closed subspaces $\{V_j:j\in \mathbb Z\}$ of $L^2(K)$ satisfying the following properties:

\parindent=0mm \vspace{.1in}
(a) $V_j \subset V_{j+1}\; \text{for all}\; j \in \mathbb Z;$

\parindent=0mm \vspace{.1in}
(b) $\bigcup_{j\in \mathbb Z}V_j\;\text{is dense in}\;L^2(K);$

\parindent=0mm \vspace{.05in}
(c) $\bigcap_{j\in \mathbb Z}V_j=\{0\};$

\parindent=0mm \vspace{.05in}
(d) $f(\cdot) \in V_j\; \text{if and only if}\;f({\mathfrak p}^{-1}\cdot) \in V_{j+1}\; \text{for all}\; j \in \mathbb Z;$

\parindent=0mm \vspace{.1in}
(e) there is a function $\varphi \in V_0$, called the {\it{scaling function}}, such that $\left\{\varphi\big(\cdot-u(k)\big): k\in \mathbb N_0\right\}$ forms a Riesz basis for $V_0$.

\parindent=8mm \vspace{.2in}
Since $\varphi \in V_0 \subset V_1\;\mbox{and}\;\{\varphi_{1,k}:k\in \mathbb N_0\} $ is a Riesz basis of $V_1$, we have

$$\varphi(x)=\sqrt{q}\sum_{k\in \mathbb N_0}a_k\, \varphi\big({\mathfrak p}^{-1}x-u(k)\big), \eqno(2.6)$$

\parindent=0mm \vspace{.1in}
where $a_k=\left<\varphi,\varphi_{1,k}\right>\;\text{and}\;\{a_k\}_{k\in \mathbb N_0}\in \ell^2(\mathbb N_0)$. Taking Fourier transform of equation (2.6), we get
\begin{align*}
\hat \varphi(x)&=\dfrac{1}{\sqrt{q}}\displaystyle\sum_{k\in \mathbb N_0}a_k \;\overline{\chi_k({\mathfrak p}\xi)}\hat \varphi({\mathfrak p}\xi)\\\\
&=m_0(\mathfrak p\xi)\,\hat \varphi(\mathfrak p\xi),
\tag{2.7}
\end{align*}
where $m_0(\xi)=\dfrac{1}{\sqrt{q}}\displaystyle\sum_{k\in \mathbb N_0}a_k\, \overline{\chi_k(\xi)}$.

\parindent=8mm \vspace{.2in}
Let $W_j,\,j\in \mathbb Z$ be the direct complementary subspace of $V_j$ in $V_{j+1}$. Assume that there exists $q-1$ functions $\{\psi_1,\psi_2,\dots,\psi_{q-1}\}$ in $L^2(K)$ such that their translates and dialations form Riesz bases of  $W_j$, i.e.,

 $$W_j=\overline{\text{span}}\left\{q^{j/2}\,\psi_\ell\big({\mathfrak p}^{-j}\cdot-u(k)\big)\right\},\,  j\in \mathbb Z,\;k\in \mathbb N_0,\;1 \le \ell \le q-1.\eqno(2.8)$$

\parindent=0mm \vspace{.1in}
Since $\psi_\ell \in W_0 \subset V_1,\;1 \le \ell \le q-1$, there exist a sequence $\{a_k^\ell\}\in \ell^2(\mathbb N_0)$ such that

$$\psi_\ell(x)=\sum_{k\in \mathbb N_0}a_k^\ell\, q^{1/2}\,\varphi\big({\mathfrak p}^{-1}x-u(k)\big),\;1 \le \ell \le q-1. \eqno(2.9)$$

\parindent=0mm \vspace{.1in}
Equation (2.9) can be written in the frequency domain as
\begin{align*}
\hat \psi_\ell(x)&=\dfrac{1}{\sqrt{q}}\displaystyle\sum_{k\in \mathbb N_0}a_k^\ell \,\overline{\chi_k({\mathfrak p}\xi)}\hat \varphi({\mathfrak p}\xi)\\\\
&=m_\ell({\mathfrak p}\xi)\,\hat \varphi({\mathfrak p}\xi),
\tag{2.10}
\end{align*}
where $m_\ell(\xi)=\dfrac{1}{\sqrt{q}}\displaystyle\sum_{k\in \mathbb N_0}a_k^\ell\,\overline{\chi_k(\xi)},\;1 \le \ell \le q-1 $.

\parindent=0mm \vspace{.2in}
{\bf{Definition 2.3.}} Let $f,\,\hat f \in L^2(K)$ be given. We say that they are {\it{biorthogonal}} if
$$\left\langle f(\cdot),\hat f\big(\cdot-u(k)\big)\right\rangle=\delta_{0,k},\eqno(2.11)$$
where
$\delta_{0,k}$ is the Kronecker's delta function.

\parindent=0mm \vspace{.2in}
If $\varphi(\cdot),\hat \varphi(\cdot) \in L^2(K)$ are a pair of biorthogonal scaling functions, then we have
$$\Big\langle\varphi(\cdot),\hat \varphi\big(\cdot-u(k)\big)\Big\rangle=\delta_{0,k}, \quad k\in\mathbb N_0.\eqno(2.12)$$

\parindent=8mm \vspace{.0in}
Further, we say that $\psi_\ell(\cdot),\hat \psi_\ell(\cdot)\in L^2(K),1 \le \ell \le q-1$ are a pair of biorthogonal wavelets associated with a pair of biorthogonal scaling functions $\varphi(\cdot),\hat \varphi(\cdot) \in L^2(K)$ if, the set $\{\psi_\ell\big(\cdot-u(k)\big): k\in \mathbb N_0, 1 \le \ell \le q-1\}$ forms a Riesz basis of $W_0$, and
\begin{align*}
\left\langle\varphi(\cdot),\tilde \psi_\ell\big(\cdot-u(k)\big)\right\rangle &=0,\quad k\in \mathbb N_0,\;1 \le \ell \le q-1,\tag{2.13}\\
\Big\langle\tilde\varphi(\cdot), \psi_\ell\big(\cdot-u(k)\big)\Big\rangle&=0,\quad  k\in \mathbb N_0,\;1 \le \ell \le q-1,\tag{2.14}\\
\left\langle \psi_\ell(\cdot),\tilde \psi_{\ell^\prime}\big(\cdot-u(k)\big)\right\rangle&=\delta_{\ell,{\ell}^\prime}\,\delta_{0,k},\quad  k\in \mathbb N_0,\;1 \le \ell,{\ell}^\prime \le q-1.\tag{2.15}
\end{align*}

\parindent=0mm \vspace{.1in}
For $\ell=1,2,\dots,q-1,$ we have

 $$W_j^\ell=\overline{\text{span}}\left\{q^{j/2}\psi_\ell\big({\mathfrak p}^{-j}\cdot-u(k)\big)\right\},\quad  j\in \mathbb Z,\;k\in \mathbb N_0.\eqno(2.16)$$

\parindent=0mm \vspace{.1in}
Using the definition of $W_j$ and identities (2.13)-(2.15), we have the following result:

\parindent=0mm \vspace{.2in}
{\bf{Proposition 2.4 [1].}} {\it{If $\psi_\ell(\cdot),\hat \psi_\ell(\cdot)\in L^2(K),\;1 \le \ell \le q-1$ are a pair of biorthogonal wavelets associated with a pair of biorthogonal scaling functions $\varphi(\cdot),\hat \varphi(\cdot) \in L^2(K)$, then }}
$$L^2(K)=\bigoplus_{j\in \mathbb Z}W_j=\bigoplus_{j\in \mathbb Z}\bigoplus_{\ell=1}^{q-1}W_j.\eqno(2.17)$$

\parindent=8mm \vspace{.1in}
In the biorthogonal setting, the refinement equation and wavelet equation are much similar to the equations (2.6) and (2.9)

$$\tilde \varphi(x)=\sqrt{q}\sum_{k\in \mathbb N_0}\tilde a_k\, \tilde\varphi\big({\mathfrak p}^{-1}x-u(k)\big),~\qquad\qquad\qquad\eqno(2.18)$$

\parindent=0mm \vspace{.0in}
and
$$\tilde\psi_\ell(x)=\sqrt{q}\sum_{k\in \mathbb N_0}\tilde a_k^\ell\,\tilde\varphi\big({\mathfrak p}^{-1}x-u(k)\big),\quad 1 \le \ell \le q-1.\eqno(2.19)$$

\parindent=0mm \vspace{.1in}
Taking Fourier transform of equations (2.18) and (2.19), we obtain

\begin{align*}
\hat{ \tilde\varphi}(x)&=\dfrac{1}{\sqrt{q}}\displaystyle\sum_{k\in \mathbb N_0}\tilde a_k \,\overline{\chi_k({\mathfrak p}\xi)}\,\hat{ \tilde\varphi}({\mathfrak p}\xi)\\\\
&=\tilde m_0({\mathfrak p}  \xi)\hat {\tilde\varphi}({\mathfrak p}\xi),\tag{2.20}
\end{align*}

where $\tilde m_0(\xi)=\dfrac{1}{\sqrt{q}}\displaystyle\sum_{k\in \mathbb N_0}\tilde a_k\, \overline{\chi_k(\xi)}$, and
\begin{align*}
\hat{\tilde \psi}_\ell(x)&=\dfrac{1}{\sqrt{q}}\displaystyle\sum_{k\in \mathbb N_0}\tilde a_k^\ell\, \overline{\chi_k({\mathfrak p}\xi)}\,\hat {\tilde\varphi}({\mathfrak p}\xi)\\\\\;
&=\tilde m_\ell({\mathfrak p} \xi)\hat{\tilde \varphi}({\mathfrak p}\xi),\tag{2.21}
\end{align*}

where $\tilde m_\ell(\xi)=\dfrac{1}{\sqrt{q}}\displaystyle\sum_{k\in \mathbb N_0}\tilde a_k^\ell \,\overline{\chi_k(\xi)},\quad 1 \le \ell \le q-1 $.

\parindent=8mm \vspace{.2in}

Moreover, it is proved in [1] that if $\varphi(\cdot),\hat \varphi(\cdot) \in L^2(K)$ are a pair of biorthogonal scaling functions associated with the MRA, then the system of functions $\{\varphi\big(\cdot-u(k)\big):k\in \mathbb N_0\}$ is biorthogonal to $\{\tilde \varphi\big(\cdot-u(k)\big):k\in \mathbb N_0\}$ if and only if
$$\displaystyle\sum_{k\in \mathbb N_0}\hat \varphi\big(\xi+u(k)\big)\overline{\hat{\tilde \varphi}\big(\xi+u(k)\big)}=1\quad a.e. \eqno(2.22)$$

\parindent=0mm \vspace{.1in}
{\bf{3. Biorthogonal Wavelet Packets on Local Fields}}

\parindent=8mm \vspace{.2in}
 For $n = 0,1,\dots$, the {\it{basic  wavelet packets}}  associated with a scaling function $\varphi(\cdot)$  on a local field $K$ of positive characteristic are defined recursively by
$$\omega_n(x)=\omega_{qr+s}(x)= \sqrt{q}\sum_{k\in \mathbb N_0} a_k^s\; \omega_r\big({\mathfrak p}^{-1}x-u(k)\big),\;\; 0 \le s \le q-1 \eqno(3.1)$$
where $r \in \mathbb N_0$ is the unique element such that $n=qr+s,\, 0\le s \le q-1$ holds (see [2]).

\parindent=8mm \vspace{.2in}
Similar to the orthogonal wavelet packets, the {\it{biorthogonal wavelet packets}} associated with $\tilde\varphi(\cdot)$ are given by

$$\tilde\omega_n(x)=\tilde\omega_{qr+s}(x)= \sqrt{q}\sum_{k\in \mathbb N_0}\tilde a_k^s\;\tilde \omega_r\big({\mathfrak p}^{-1}x-u(k)\big),\;\; 0 \le s \le q-1. \eqno(3.2)$$

\parindent=0mm \vspace{.1in}
Clearly, for $r=0$ and $1 \le s \le q-1$, we have
$$\omega_0(\cdot)=\varphi(\cdot),\quad\tilde\omega_0(\cdot)=\tilde\varphi(\cdot),\quad\omega_s(\cdot)=\psi_\ell(\cdot),\quad\tilde\omega_s(\cdot)=\tilde\psi_\ell(\cdot), \quad 1\le \ell \le q-1.$$

\parindent=0mm \vspace{.1in}
Also, the Fourier transform of (3.1) and (3.2) gives

$$\hat \omega_{qr+s}(\xi)=m_s({\mathfrak p}\xi)\,\hat\omega_r({\mathfrak p}\xi), \eqno(3.3)$$
and
$$\hat {\tilde \omega}_{qr+s}(\xi)=\tilde m_s({\mathfrak p}\xi)\,\hat {\tilde \omega}_r({\mathfrak p}\xi). \eqno(3.4)$$

\parindent=8mm \vspace{.1in}
We are now in a position to discuss the biorthogonality properties for these wavelet packets by means of the Fourier transform.

\parindent=0mm \vspace{.2in}
{\bf{Lemma 3.1.}} {\it{Assume that $\omega_s(x),\,\tilde \omega_s(x) \in L^2(K)$ are a pair of biorthogonal wavelets associated with a pair of biorthogonal scaling functions $\omega_0(x),\,\tilde \omega_0(x)$. Then we have}}
$$\sum_{\ell=0}^{q-1}m_r\big({\mathfrak p}\xi+{\mathfrak p} u(\ell)\big)\overline{\tilde m_s\big({\mathfrak p}\xi+{\mathfrak p} u(\ell)\big)}= \delta_{r, s},\quad 0 \le r, s \le q-1.\eqno(3.5)$$

\parindent=0mm \vspace{.1in}
{\bf{Proof.}} For given $0 \le r, s \le q-1$, we have

\parindent=0mm \vspace{.2in}

$\begin{array}{rcl}
\delta_{r, s}&=&\displaystyle\sum_{k\in \mathbb N_0}\omega_r\big(\xi+ u(k)\big)\overline{\tilde \omega_r\big(\xi+ u(k)\big)}\\\\
&=&\displaystyle\sum_{k\in \mathbb N_0}m_r\big({\mathfrak p}\xi+{\mathfrak p} u(k)\big)\hat \omega_0\big({\mathfrak p}\xi+{\mathfrak p} u(k)\big)\overline{\hat{\tilde \omega}_0\big({\mathfrak p}\xi+{\mathfrak p} u(k)\big)} \;\overline{ \tilde m_s\big({\mathfrak p}\xi+{\mathfrak p} u(k)\big)}\\\\
&=&\displaystyle\sum_{\ell=0}^{q-1}\displaystyle\sum_{k\in \mathbb N_0}m_r\big({\mathfrak p}\xi+{\mathfrak p} u(qk+\ell)\big)\hat \omega_0\big({\mathfrak p}\xi+{\mathfrak p} u(qk+\ell)\big)\overline{\hat{\tilde \omega}_0\big({\mathfrak p}\xi+{\mathfrak p} u(qk+\ell)\big)}\\
&&\qquad\qquad\qquad\qquad\qquad\qquad\qquad\qquad\qquad\qquad\times\quad \overline{ \tilde m_s\big({\mathfrak p}\xi+{\mathfrak p} u(qk+\ell)\big)}\\\\
&=&\displaystyle\sum_{\ell=0}^{q-1}m_r\big({\mathfrak p}\xi+{\mathfrak p} u(\ell)\big)\overline{m_s\big({\mathfrak p}\xi+{\mathfrak p} u(\ell)\big)}\\
&&\qquad\qquad\qquad\qquad \times \left\{\displaystyle\sum_{k\in \mathbb N_0}\hat \omega_0\big({\mathfrak p}\xi+{\mathfrak p} u(qk+\ell)\big) \overline{\hat{\tilde \omega}_0\big({\mathfrak p}\xi+{\mathfrak p} u(qk+\ell)\big)  }\right\}\\\\
&=&\displaystyle\sum_{\ell=0}^{q-1}m_r\big({\mathfrak p}\xi+{\mathfrak p} u(\ell)\big)\overline{m_s\big({\mathfrak p}\xi+{\mathfrak p} u(\ell)\big)}.
\end{array}$

\parindent=0mm \vspace{.2in}
{\bf{Theorem 3.2.}}{\it{ If $\left\{\omega_n(x):n\in \mathbb N_0\right\}$ and $\left\{\tilde\omega_n(x):n\in \mathbb N_0\right\}$ are wavelet packets associated with a pair of biorthogonal scaling functions $\omega_0(x)$ and $\tilde \omega_0(x)$, respectively. Then, we have}}

$$\Big\langle\omega_n(\cdot),\tilde\omega_n\big(\cdot-u(k)\big)\Big\rangle=\delta_{0,k},\quad  k\in \mathbb Z,\,n \in \mathbb N_0.\eqno(3.6)$$

\parindent=0mm \vspace{.1in}
{\bf{Proof.}} We will prove this result by using induction on $n$. It follows from (2.12) and (2.15) that the claim is true for $n=0$ and $n=1,2,\dots,q-1$. Assume (3.6) holds for $n<t$, where $t\in \mathbb N$. Then, we prove the result (3.6) for $n=t$. Let $n=qr+s$, where $r\in \mathbb N_0,\,0 \le s \le q-1$, and $r<n$. Therefore, by the inductive assumption, we have
$$\Big\langle\omega_r(\cdot),\tilde\omega_r\big(\cdot-u(k)\big)\Big\rangle=\delta_{0,k}\;\Longleftrightarrow\; \sum_{k\in \mathbb N_0}\omega_r\big(\xi+ u(k)\big)\,\overline{\tilde \omega_r\big(\xi+ u(k)\big)}=1.$$

Using Lemmas 2.5, 3.1 and, equations (3.3) and (3.4), we obtain

\parindent=0mm \vspace{.2in}
$ \Big\langle\omega_n(\cdot),\tilde\omega_n\big(\cdot-u(k)\big)\Big\rangle\\\\
\begin{array}{rcl}
&&=\quad \left\langle \hat \omega_n(\cdot),\hat{\tilde\omega}_n\big(\cdot-u(k)\big)\right\rangle\\\\
 &&=\quad \displaystyle\int_K \hat \omega_{qr+s}(\xi)\;\overline{\hat{\tilde \omega}_{qr+s}(\xi)}\,\chi_k(\xi)d \xi\\\\
&&=\quad \displaystyle\int_K m_s({\mathfrak  p} \xi)\hat \omega_r({\mathfrak  p }\xi)\;\overline{\tilde m_s({\mathfrak  p }\xi)}\;\overline{\hat{\tilde \omega}_r({\mathfrak  p }\xi)}\,\chi_k(\xi)d \xi\\\\
&&=\quad \displaystyle\int_{\mathfrak D} \displaystyle\sum_{k\in \mathbb N_0}m_s\big({\mathfrak  p }\xi + {\mathfrak  p }u(k)\big)\,\hat \omega_r({\mathfrak  p }\xi + {\mathfrak  p }u(k))\\\
&& \qquad\qquad\qquad\qquad\qquad\qquad\times~\overline{\tilde m_s\big({\mathfrak  p }\xi + {\mathfrak  p }u(k)\big)}\,\overline{\hat{\tilde \omega}_r\big({\mathfrak  p }\xi + {\mathfrak  p }u(k)\big)}\,\chi_k(\xi)d \xi\\\\
&&=\quad \displaystyle\int_{\mathfrak D} \displaystyle\sum_{\ell=0}^{q-1}\displaystyle\sum_{k\in \mathbb N_0}m_s\big({\mathfrak  p }\xi + {\mathfrak  p }u(qk+\ell)\big)\,\hat \omega_r\big({\mathfrak  p }\xi + {\mathfrak  p }u(qk+\ell)\big)\\\
&&\qquad\qquad\qquad\qquad\qquad   \times~ \overline{\tilde m_s\big({\mathfrak  p }\xi + {\mathfrak  p }u(qk+\ell)\big)}\,\overline{\hat{\tilde \omega}_r\big({\mathfrak  p }\xi + {\mathfrak  p }u(qk+\ell)\big)}\,\chi_k(\xi)d \xi\\\\
&&=\quad \displaystyle\int_{\mathfrak D} \displaystyle\sum_{\ell=0}^{q-1}m_s\big({\mathfrak  p }\xi + {\mathfrak  p }u(\ell)\big)\,\overline{\tilde m_s\big({\mathfrak  p }\xi + {\mathfrak  p }u(\ell)\big)}\\\
&&\qquad\qquad\qquad\qquad\times~ \left\{\displaystyle\sum_{k\in \mathbb N_0}\,\hat \omega_r\big({\mathfrak  p }\xi + {\mathfrak  p }u(qk+\ell)\big)\,\overline{\hat{\tilde \omega}_r\big({\mathfrak  p }\xi + {\mathfrak  p }u(qk+\ell)\big)}\,\chi_k(\xi)d \xi\right\}\\\\
&&=\quad \displaystyle\int_{\mathfrak D} \displaystyle\sum_{\ell=0}^{q-1}m_s\big({\mathfrak  p }\xi + {\mathfrak  p }u(\ell)\big)\,\overline{\tilde m_s\big({\mathfrak  p }\xi + {\mathfrak  p }u(\ell)\big)}\,\chi_k(\xi)d \xi\\\\
&&=\quad \delta_{0,k}.
\end{array}$

\parindent=0mm \vspace{.3in}
{\bf{Theorem 3.3.}}{\it{ Suppose $\left\{\omega_n(x):n\in \mathbb N_0\right\}$ and $\left\{\tilde\omega_n(x):n\in \mathbb N_0\right\}$ are the biorthogonal wavelet packets associated with a pair of biorthogonal scaling functions $\omega_0(x)$ and $\tilde \omega_0(x)$, respectively. Then, we have}}

$$\Big\langle\omega_{qr+s_1}(\cdot),\tilde\omega_{qr+s_2}\big(\cdot-u(k)\big)\Big\rangle=\delta_{0,k}\delta_{s_1,s_2} ,\;0 \le \,s_1,s_2 \,\le q-1,\; r,k\in \mathbb N_0.\eqno(3.7)$$

\parindent=0mm \vspace{.2in}
{\bf{Proof.}} By Lemma 2.5, we have\\

$\Big\langle\omega_{qr+s_1}(\cdot),\tilde\omega_{qr+s_2}\big(\cdot-u(k)\big)\Big\rangle$ \\

$\begin{array}{lrl}
&&= \left\langle\hat \omega_{qr+s_1}(\cdot),\hat{\tilde\omega}_{qr+s_2}\big(\cdot-u(k)\big)\right\rangle\\\\
&&=\displaystyle\int_K \hat \omega_{qr+s_1}(\xi)\,\overline{\hat{\tilde \omega}_{qr+s_2}(\xi)}\,\chi_k(\xi)d \xi\\\\
&&=\quad \displaystyle\int_K m_{s_1}({\mathfrak  p }\xi)\,\hat \omega_r({\mathfrak  p }\xi)\,\overline{\tilde m_{s_2}({\mathfrak  p }\xi)}\, \overline{\hat{\tilde \omega}_r({\mathfrak  p }\xi)}\,\chi_k(\xi)d \xi\\\\
&&=\quad \displaystyle\int_{\mathfrak D} \displaystyle\sum_{k\in \mathbb N_0}m_{s_1}({\mathfrak  p }\xi + {\mathfrak  p }u(k))\,\hat \omega_r({\mathfrak  p }\xi + {\mathfrak  p }u(k))\\
&&\qquad\qquad\qquad\qquad\qquad \qquad\quad\times ~\overline{\tilde m_{s_2}({\mathfrak  p }\xi + {\mathfrak  p }u(k))}\,\, \overline{\hat{\tilde \omega}_r({\mathfrak  p }\xi + {\mathfrak  p }u(k))}\,\chi_k(\xi)d \xi\\\\
&&=\quad \displaystyle\int_{\mathfrak D} \displaystyle\sum_{\ell=0}^{q-1}\displaystyle\sum_{k\in \mathbb N_0}m_{s_1}\big({\mathfrak  p }\xi + {\mathfrak  p }u(qk+\ell)\big)\hat \omega_r\big({\mathfrak  p }\xi + {\mathfrak  p }u(qk+\ell)\big)\\
&&\qquad\qquad\qquad\qquad\qquad\times~ \overline{\tilde m_{s_2}\big({\mathfrak  p }\xi + {\mathfrak  p }u(qk+\ell)\big)} \,\,\overline{\hat{\tilde \omega}_r\big({\mathfrak  p }\xi + {\mathfrak  p }u(qk+\ell)\big)}\chi_k(\xi)d \xi\\\\
&&=\quad \displaystyle\int_{\mathfrak D} \displaystyle\sum_{\ell=0}^{q-1}m_{s_1}\big({\mathfrak  p }\xi + {\mathfrak  p }u(\ell)\big)\,\overline{\tilde m_{s_2}\big({\mathfrak  p }\xi + {\mathfrak  p }u(\ell)\big)}\\
&&\quad \quad  \quad \quad \quad \quad \quad \quad \quad \quad \times \left\{\displaystyle\sum_{k\in \mathbb N_0}\hat \omega_r\big({\mathfrak  p }\xi + {\mathfrak  p }u(qk+\ell)\big)\,\,\overline{\hat{\tilde \omega}_r\big({\mathfrak  p }\xi + {\mathfrak  p }u(qk+\ell)\big)}\,\chi_k(\xi)d \xi\right\}\\\\
&&=\quad \displaystyle\int_{\mathfrak D} \displaystyle\sum_{\ell=0}^{q-1}m_{s_1}\big({\mathfrak  p }\xi + {\mathfrak  p }u(\ell)\big)\,\overline{\tilde m_{s_2}\big({\mathfrak  p }\xi + {\mathfrak  p }u(\ell)\big)}\,\chi_k(\xi)d \xi\\\\
&&=\quad \delta_{0,k}\delta_{s_1,s_2}.
\end{array}$

\parindent=0mm \vspace{.2in}
{\bf{Theorem 3.4.}}{\it{ Suppose $\left\{\omega_n(x):n\in \mathbb N_0\right\}$ and $\left\{\tilde\omega_n(x):n\in \mathbb N_0\right\}$ are wavelet packets with respect to a pair of biorthogonal scaling functions $\omega_0(x)$ and $\tilde \omega_0(x)$, respectively. Then, we have}}
$$\Big\langle\omega_\ell(\cdot),\tilde\omega_n\big(\cdot-u(k)\big)\Big\rangle=\delta_{\ell,n}\,\delta_{0,k},\quad \ell,n, k\in \mathbb N_0.\eqno(3.8)$$

\parindent=0mm \vspace{.1in}
{\bf{Proof.}} For $\ell=n$, the result (3.8) follows by Theorem 3.2. When $\ell \neq n$, and $0 \le \ell,n \le q-1$, the result (3.8) can be established from Theorem 3.3. Assume $\ell$ is not equal to $n$ and at least one of $\ell,n$ does not belong to $\left\{1,2,\dots,q-1\right\}$, then we can write $\ell,n$ as $\ell=qr_1+s_1,\;n=qu_1+v_1,\,r_1,u_1 \in \mathbb N_0,\;s_1,v_1\in \left\{0,1,2,\dots,q-1\right\}$.

\parindent=0mm \vspace{.2in}
{\bf{Case 1:}} If $r_1 = u_1$, then $s_1 \neq v_1$. Therefore, (3.8) follows by virtue of the properties (3.3)-(3.5) and Lemma 2.5 i.e.,

\parindent=0mm \vspace{.2in}
$\Big\langle\omega_\ell(\cdot),\tilde\omega_n\big(\cdot-u(k)\big)\Big\rangle\\\\
\begin{array}{lcl}
&=&\Big\langle\omega_{qr_1+s_1}(\cdot),\tilde\omega_{qu_1+v_1}\big(\cdot-u(k)\big)\Big\rangle\\\\
&=&\left\langle\hat \omega_{qr_1+s_1}(\cdot),\hat{\tilde\omega}_{qu_1+v_1}\big(\cdot-u(k)\big)\right\rangle\\\\
&=& \displaystyle\int_K \hat \omega_{qr_1+s_1}(\xi)\,\overline{\hat{\tilde \omega}_{qu_1+v_1}(\xi)}\chi_k(\xi)d \xi\\\\
&=& \displaystyle\int_K m_{s_1}({\mathfrak  p }\xi)\,\hat \omega_{r_1}({\mathfrak  p }\xi)\,\overline{\tilde m_{v_1}({\mathfrak  p }\xi)}\, \overline{\hat{\tilde \omega}_{u_1}({\mathfrak  p }\xi)}\,\chi_k(\xi)d \xi\\\\
&=& \displaystyle\int_{\mathfrak D} \displaystyle\sum_{k\in \mathbb N_0}m_{s_1}\big({\mathfrak  p }\xi + {\mathfrak  p }u(k)\big)\,\hat \omega_{r_1}({\mathfrak  p }\xi + {\mathfrak  p }u(k))\\
&&~~~~~~~~~~~~~~~~~~~~~~~~~~\times \overline{\tilde m_{v_1}\big({\mathfrak  p }\xi + {\mathfrak  p }u(k)\big)} \,\overline{\hat{\tilde \omega}_{u_1}\big({\mathfrak  p }\xi + {\mathfrak  p }u(k)\big)}\chi_k(\xi)d \xi\\\\
&=& \displaystyle\int_{\mathfrak D} \displaystyle\sum_{\ell=0}^{q-1}\displaystyle\sum_{k\in \mathbb N_0}m_{s_1}\big({\mathfrak  p }\xi + {\mathfrak  p }u(qk+\ell)\big)\,\hat \omega_{r_1}\big({\mathfrak  p }\xi + {\mathfrak  p }u(qk+\ell)\big)\\
&& ~~~~~~~~~~~~~~~~~~~~~~~~~~\times \overline{\tilde m_{v_1}\big({\mathfrak  p }\xi + {\mathfrak  p }u(qk+\ell)\big)}\, \overline{\hat{\tilde \omega}_{u_1}\big({\mathfrak  p }\xi + {\mathfrak  p }u(qk+\ell)\big)}\,\chi_k(\xi)d \xi\\\\
&=& \displaystyle\int_{\mathfrak D} \displaystyle\sum_{\ell=0}^{q-1}m_{s_1}\big({\mathfrak  p }\xi+ {\mathfrak p} u(\ell)+{\mathfrak p} u(\ell)\big)\,\overline{\tilde m_{v_1}\big({\mathfrak  p }\xi+ {\mathfrak p} u(\ell)\big)}\\
&&~~~~~~~~~~~~~~~~~~~~~~~~~~\times \left\{\displaystyle\sum_{k\in \mathbb N_0}\hat \omega_{r_1}\big({\mathfrak  p }\xi + {\mathfrak  p }u(qk+\ell)\big)\,\overline{\hat{\tilde \omega}_{u_1}\big({\mathfrak  p }\xi + {\mathfrak  p }u(qk+\ell)\big)}\chi_k(\xi)d \xi\right\}\\\\
&=& \displaystyle\int_{\mathfrak D} \displaystyle\sum_{\ell=0}^{q-1}m_{s_1}\big({\mathfrak  p }\xi+ {\mathfrak p} u(\ell)\big)\,\overline{\tilde m_{v_1}\big({\mathfrak  p }\xi+ {\mathfrak p} u(\ell)\big)}\,\chi_k(\xi)d \xi\\\\
&=&\delta_{0,k}.
\end{array}$

\parindent=0mm \vspace{.2in}
{\bf{Case 2:}} If $r_1 \neq u_1$, then $r_1=pr_2+s_2,\, u_1=pu_2+v_2$, where $r_2,u_2 \in \mathbb N_0,$ and $ s_2,v_2\in \left\{0,1,\dots,q-1\right\}$. If $r_2 = u_2$, then $s_2 \neq v_2$. Similar to Case 1, (3.8) can be established. When  $r_2 \neq u_2$, we order $r_2=pr_3+s_3,\, u_2=pu_3+v_3,$ where $r_3,u_3 \in \mathbb N_0,$ and $s_3,v_3\in \left\{0,1,\dots,q-1\right\}$. Thus, after taking finite steps (denoted by $h$), we obtain $r_h,u_h \in \mathbb N_0$ and $s_h,v_h\in \left\{0,1,\dots,q-1\right\}$. If $r_h = u_h$, then $s_h \neq v_h$. Similar to Case 1, (3.8) can be established. When   $r_h \neq u_h$, it follows from equations (2.12)-(2.15) that

$$\Big\langle\omega_{r_h}(\cdot),\tilde\omega_{u_h}\big(\cdot-u(k)\big)\Big\rangle=0\;\Longleftrightarrow\; \sum_{k\in \mathbb N_0}\omega_{r_h}\big(\xi+ u(k)\big)\,\overline{\tilde \omega_{u_h}\big(\xi+ u(k)\big)}=0,\;\; \xi \in  K.$$

Also, we have\\\\
$\Big\langle\omega_r(\cdot),\tilde\omega_u\big(\cdot-u(k)\big)\Big\rangle\\\\
\begin{array}{lcl}
&=&\left\langle\hat \omega_r(\cdot),\hat{\tilde\omega}_u\big(\cdot-u(k)\big)\right\rangle\\\\
&=&\left\langle\hat\omega_{qr_1+s_1}(\cdot),\hat{\tilde\omega}_{qu_1+v_1}\big(\cdot-u(k)\big)\right\rangle\\\\
&=& \displaystyle\int_K \hat \omega_{qr_1+s_1}(\xi)\,\overline{\hat{\tilde \omega}_{qu_1+v_1}(\xi)}\chi_k(\xi)d \xi\\\\
&=& \displaystyle\int_K m_{s_1}({\mathfrak  p }\xi)\,m_{s_2}({\mathfrak  p }^2 \xi)\,\hat \omega_{r_2}({\mathfrak  p }^2 \xi)\,\overline{\tilde m_{v_1}({\mathfrak  p }\xi)}\,\overline{\tilde m_{v_2}({\mathfrak  p }^2 \xi)}\,
\overline{\hat{\tilde \omega}_{u_2}({\mathfrak  p }^2 \xi)}\,\chi_k(\xi)d \xi\\\\
&& \vdots\\\\
&=& \displaystyle\int_K\left\{ \prod_{\ell=1}^h m_{s_\ell}({\mathfrak  p }^\ell \xi)\right\}\hat \omega_{r_h}({\mathfrak  p }^h \xi)\,\left\{ \prod_{\ell=1}^h\overline{\tilde m_{v_\ell}({\mathfrak  p }^\ell \xi)}\right\}\overline{\hat{\tilde \omega}_{u_h}({\mathfrak  p }^h \xi)}\,\chi_k(\xi)d \xi\\\\
&=& \displaystyle\int_{\mathfrak D} \displaystyle\sum_{k\in \mathbb N_0}\left\{ \prod_{\ell=1}^h m_{s_\ell}\left({\mathfrak  p }^\ell\big(\xi+ u(k)\big)\right)\right\}\left\{\hat \omega_{r_h}\Big({\mathfrak  p }^h\big(\xi+ u(k)\big)\Big)\,\overline{\hat \omega_{u_h}\Big({\mathfrak  p }^h\big(\xi+ u(k)\big)\Big)}\right\}\\
&&\qquad\qquad\qquad\qquad\qquad\qquad\qquad\qquad\times \left\{\displaystyle\prod_{\ell=1}^h m_{v_\ell}\Big({\mathfrak  p }^\ell\big(\xi+ u(k)\big)\Big)\right\}\chi_k(\xi)d \xi\\
&=&0.
\end{array}$

\parindent=0mm \vspace{.3in}

 {\bf{4. Construction of Riesz Bases from Wavelet Packets }}

\parindent=8mm \vspace{.2in}
In this section, we will decompose the subspaces  $V_j,  \tilde V_j,\,W_j\, \text{and}\,\tilde W_j$ by constructing subspaces of wavelet packets. We also present a direct decomposition for $L^2(K)$.

\parindent=8mm \vspace{.2in}
For any $n\in \mathbb N_0$, define
$$E_n=\left\{f(x): f(x)=\sum_{k\in \mathbb N_0} a_k \omega_n(x- u(k)),~ \left\{ a_k\right\}_{k\in \mathbb N_0}\in l^2(\mathbb N_0)   \right\},\eqno(4.1)$$
~$$\tilde E_n=\left\{\tilde f(x): \tilde f(x)=\sum_{k\in \mathbb N_0} \tilde a_k \tilde \omega_n(x- u(k)),~ \left\{ \tilde a_k\right\}_{k\in \mathbb N_0}\in l^2(\mathbb N_0)   \right\}.\eqno(4.2)$$
Clearly $E_0=V_0$ and $E_s=W_0^s$, for any $1 \le s \le q-1$. Assume that $\Big\{m_s\big(p\xi + {\mathfrak  p } u(k)\big)\Big\}_{s,k=0}^{q-1}$ is a unitary matrix.

\parindent=0mm \vspace{.2in}
{\bf{Lemma 4.1.}} {\it{For $n \in \mathbb N_0$, the space $\Delta E_n$ can be decomposed into the direct sum of  $U_{qn+s}, 1 \le s \le q-1$, i.e.,
$$\Delta E_n=\displaystyle\bigoplus_{s =0}^{q-1} U_{qn+s}, \eqno(4.3)$$

where $\Delta$ is the dialation operator such that $\Delta f(x)=f({\mathfrak  p }^{-1}x)$, for any $f \in L^2(K)$.}}

\parindent=0mm \vspace{.2in}
{\bf{Proof.}} For $n \in \mathbb N_0$, we claim that
~~$$\Delta E_n=\left\{f(x): f(x)=\sum_{s =0}^{q-1}\sum_{k\in \mathbb N_0} a_k^s \;\omega_{qn+s}\big(x- u(k)\big),~ \left\{ a_k^s\right\}_{k\in \mathbb N_0}\in l^2(\mathbb N_0)   \right\}.\eqno(4.4)$$

\parindent=0mm \vspace{.1in}
As for any $0 \le s \le q-1$, by $(3.1)$ and $(4.1)$, $\omega_{qn+s}\big(x- u(k)\big) \in \Delta E_n$. Assume that $f(x) \in \Delta E_n$, then there exists a sequence $\left\{ b_k\right\}_{k\in \mathbb N_0}\in l^2(\mathbb N_0)$ such that
$$f(x)= \sum_{k\in \mathbb N_0} b_k \;\omega_n\big({\mathfrak  p }^{-1}x-u(k)\big). \eqno(4.5)$$

Similarly, for each $s=0,1,\dots, q-1$, there exist a sequence $\left\{a_k^s\right\}_{k\in \mathbb N_0}$ in $l^2(\mathbb N_0)$ such that

$$f(x)= \sum_{s=1}^{q-1}\sum_{k\in \mathbb N_0} a_k^s\, \omega_n\big({\mathfrak  p }^{-1}x-u(k)\big), \eqno(4.6)$$

\parindent=0mm \vspace{.1in}
provided $f(x) \in \Delta E_n$.

\parindent=8mm \vspace{.2in}
Taking Fourier transform on the both sides of (4.5) and (4.6), respectively and by using (3.3), we obtain
$$\hat f(\xi)=h({\mathfrak  p }\xi)\hat \omega_n({\mathfrak  p }\xi)=\sum_{s=1}^{q-1} g_s(\xi) m_s({\mathfrak  p }\xi)\hat \omega_n({\mathfrak  p }\xi),\eqno(4.7)$$
where $h(\xi)= \displaystyle\sum_{k\in \mathbb N_0} b_k \overline{\chi_k(\xi)},\;g_s(\xi)= \displaystyle\sum_{k\in \mathbb N_0} a_k^s\overline{\chi_k(\xi)}.$\\

\parindent=0mm \vspace{.2in}
The above result (4.7) follows if the following equality holds:
$$h({\mathfrak  p }\xi)=\sum_{s=1}^{q-1}g_s(\xi)\, m_s({\mathfrak  p }\xi). \eqno(4.8)$$
For any $\left\{ b_k\right\}_{k\in \mathbb N_0}\in l^2(\mathbb N_0)$, we will prove that there exists a sequence $\left\{ a_k^s\right\}_{k\in \mathbb N_0}\in l^2(\mathbb N_0)$ such that (4.8) is satisfied. Moreover, equation(4.8) is equal to the following identity:
$$h\big({\mathfrak  p }\xi +{\mathfrak  p } u(k)\big)=\sum_{s=1}^{q-1}g_s(\xi)\, m_s\big({\mathfrak  p }\xi+{\mathfrak  p } u(k)\big). \eqno(4.9)$$
The solvibility of (4.9) for every sequence $\left\{ b_k\right\}_{k\in \mathbb N_0}\in l^2(\mathbb N_0)$ follows from the fact that the matrix $\Big\{m_s\big(p\xi + {\mathfrak  p } u(k)\big)\Big\}_{s,k=0}^{q-1}$ is unitary. Hence, equation (4.4) follows. Further, applying Theorem 3.3, it follows that

$$\left\{\omega_{qn+s}\big({\mathfrak  p }^{-1}x-u(k)\big)\;n\in \mathbb N_0,\;\;0 \le s \le q-1,\;\;k\in \mathbb N_0 \right\}$$

\parindent=0mm \vspace{.1in}
is a Riesz bases of $\Delta E_n$.

\parindent=8mm \vspace{.2in}
Similar to (4.3), we can establish the following results:

$$\tilde E_0=\tilde V_0,\;\;\tilde E_s=\tilde W_0^s,\quad 1 \le s \le q-1,$$
and
$$\Delta \tilde E_n=\displaystyle\bigoplus_{s =0}^{q-1} \tilde U_{qn+s},\quad 1 \le s \le q-1. \qquad~~\eqno(4.10)$$

\parindent=8mm \vspace{.0in}
For $\ell \in \mathbb N$, define $\tilde \vartheta_\ell=\displaystyle\sum_{j=0}^\ell q^j \Lambda,$ where $\Lambda=\left\{0,1,2,\dots,q-1\right\},\; \vartheta_\ell=\tilde \vartheta_\ell-\tilde \vartheta_{\ell-1}$. Now, we will establish the direct decomposition of the space $L^2(K)$.

\parindent=0mm \vspace{.2in}
{\bf{Theorem 4.2.}} {\it{The family of functions $\left\{ \omega_n\big(x-u(k)\big),\;n \in \vartheta_\ell,\; k \in \mathbb N_0\right\}$ constitutes Riesz basis of $\Delta^\ell W_0$. In particular $\left\{ \omega_n\big(x-u(k)\big),\;n \in \vartheta_\ell,\; k \in \mathbb N_0\right\}$ constitutes Riesz basis of $L^2(K)$.}}

\parindent=0mm \vspace{.2in}
{\bf{Proof.}} From equation (4.3), we have

$$\Delta E_0=\displaystyle\bigoplus_{s =0}^{q-1}E_s ~ \text{i.e.,}~ \Delta  E_0=E_0\displaystyle\bigoplus_{s =1}^{q-1} E_s.$$

\parindent=0mm \vspace{.1in}
Since $E_0=V_0$ and $W_0=\bigoplus_{s =1}^{q-1}W_0^s =\bigoplus_{s =1}^{q-1} E_s$, therefore, $\Delta  E_0=V_0\bigoplus W_0$. It can be inductively
inferred from (4.3) that

$$\Delta^\ell E_0=\Delta^{\ell-1} E_0\displaystyle\bigoplus_{n \in \vartheta_\ell}E_n,\quad \ell \in \mathbb N. \eqno(4.11)$$

\parindent=0mm \vspace{.1in}
Since $V_{j+1}=V_j\bigoplus W_j,\, j\in \mathbb Z$, hence, $\Delta^\ell E_0=\delta^{\ell-1} E_0\bigoplus_{n \in \vartheta}W_0,\,\ell \in \mathbb N$. Now, it follows from (4.3) and Proposition 2.4 that $\Delta^\ell W_0=\bigoplus_{n \in \vartheta}E_n,$
and

$$L^2(K)=V_0\bigoplus\left(\displaystyle\bigoplus_{\ell \geq 0}\Delta^\ell W_0\right)=E_0\bigoplus\left(\displaystyle\bigoplus_{\ell \geq 0}\left(\displaystyle\bigoplus_{n \in \vartheta_\ell}E_n\right)\right)=\displaystyle\bigoplus_{n \in \mathbb N_0}E_n. \eqno(4.12)$$

\parindent=0mm \vspace{.1in}
In view of Theorem 3.3, the family of functions $\left\{ \omega_n\big(x-u(k)\big),n \in \vartheta_\ell, k \in \mathbb N_0\right\}$ is a  Riesz basis of $\Delta^\ell W_0$. Thus, according to (4.12), the family $\left\{ \omega_n\big(x-u(k)\big),n \in \vartheta_\ell, k \in \mathbb N_0\right\}$ forms a Riesz basis of $L^2(K)$.

\parindent=0mm \vspace{.2in}
{\bf{Corollary 4.3.}} {\it{For every $\ell \in \mathbb N$, the family of functions $\left\{\tilde \omega_n\big(x-u(k)\big),n \in \vartheta_\ell, k \in \mathbb N_0\right\}$ forms a Riesz basis of $\tilde\Delta^l W_0$.}}

\parindent=0mm \vspace{.2in}
{\bf{Corollary 4.4.}}{\it{ For every $\ell \in \mathbb N$, the family of functions $\left\{ \omega_n\big({\mathfrak  p }^{-j}x-u(k)\big),n \in \vartheta_\ell, k \in \mathbb N_0\right\}$ forms a Riesz basis of $L^2(K)$.}}

\parindent=0mm \vspace{.3in}

{\bf{5. Decomposition and Reconstruction Algorithms}}

\parindent=8mm \vspace{.2in}
We begin this section with the decomposition formulae for the biorthogonal wavelet packects on local fields of positive characteristic followed by an algorithm.

\parindent=0mm \vspace{.2in}
{\bf Theorem 5.1.} {\it Let $\{\omega_n\}$ and $\{\tilde \omega_n\}$ be the biorthogonal wavelet packets defined by (3.1) and (3.2), respectively. Then for all $k \in \mathbb N_0$, we have the following decomposition formulae:
$$\omega_n\big({\mathfrak p}^{-1}x-u(k)\big)=\dfrac{1}{\sqrt{q}}\sum_{\nu=1}^{q-1}\sum_{\mu \in \mathbb N_0}\tilde a_{k-q\mu}^\nu \,\omega_{qn+\nu}\big(x-u(\mu)\big), \eqno(5.1)$$
and
$$\tilde \omega_n\big({\mathfrak p}^{-1}x-u(k)\big)=\dfrac{1}{\sqrt{q}}\sum_{\nu=1}^{q-1}\sum_{\mu \in \mathbb N_0} a_{k-q\mu}^\nu \,\tilde\omega_{qn+\nu}\big(x-u(\mu)\big). \eqno(5.2)$$ }

\parindent=0mm \vspace{.1in}
{\bf Proof.} We will prove only (5.1). The second formula (5.2) being the dual of (5.1) will follow. Thus using equation (3.1), we have\\
$$\begin{array}{lll}
\dfrac{1}{\sqrt{q}}\displaystyle\sum_{\nu=0}^{q-1}\displaystyle\sum_{\mu \in \mathbb N_0}\tilde a_{k-q\mu}^\nu \omega_{qn+\nu}\big(x-u(\mu)\big)\\\\
\quad \quad =\quad \dfrac{1}{\sqrt{q}}\displaystyle\sum_{\nu=0}^{q-1}\displaystyle\sum_{\mu \in \mathbb N_0}\tilde a_{k-q\mu}^\nu\,q^{1/2} \displaystyle\sum_{r \in \mathbb N_0} a_r^\nu \omega_n\left({\mathfrak p}^{-1}\big(x-u(\mu)\big)-u(r)\right)&&\\\\
\quad \quad =\quad \displaystyle\sum_{\nu=0}^{q-1}\displaystyle\sum_{\mu \in \mathbb N_0}\tilde a_{k-q\mu}^\nu \displaystyle\sum_{r \in \mathbb N_0} a_r^\nu \omega_n\big({\mathfrak p}^{-1}x-u(q\mu-r)\big)&&\\\\
\quad \quad =\quad \displaystyle\sum_{\nu=0}^{q-1}\displaystyle\sum_{t \in \mathbb N_0}\omega_n\big({\mathfrak p}^{-1}x-u(t)\big) \displaystyle\sum_{\mu \in \mathbb N_0}\tilde a_{k-q\mu}^\nu a_{t-q\mu}^\nu &&\\\\
\quad \quad =\quad \omega_n\big({\mathfrak p}^{-1}x-u(k)\big).
\end{array}$$

\parindent=0mm \vspace{.2in}

This completes the proof of the Theorem.

\parindent=8mm \vspace{.2in}

Given a level $J$ and consider
$$f\approx f_J= \sum_{k \in \mathbb N_0}c_k^J \omega_0\big({\mathfrak p}^{-J}x-u(k)\big),$$
where $\{c_k^J\} \in l^2(\mathbb N_0)$. Using the fact

$$V_J=W_{J-1}\oplus V_{J-1}=\cdots=W_{J-1}\oplus W_{J-2}\oplus \cdots W_{J-M}\oplus V_{J-M},$$

\parindent=0mm \vspace{.1in}
one obtains
$$f_J=g_{J-1}+g_{J-2}+\cdots+g_{J-M}+ f_{J-M},$$

\parindent=0mm \vspace{.1in}
where $f_{J-M} \in V_{J-M}$ and $g_j \in W_j,\; j= J-M, \dots, J-1.$

\parindent=0mm \vspace{.2in}
Furthermore, by using Theorem 5.1, $g_j \in W_j, j= J-M, \dots, J-1$ can be further decomposed. To do this, let

$$f_j(x)= \sum_{k \in \mathbb N_0}c_k^j \,\omega_0 \big({\mathfrak p}^{-j}x-u(k)\big), \eqno (5.3)$$
and
$$g_j(x)= \sum_{\nu=1}^{q-1}\sum_{k \in \mathbb N_0}d_k^{\nu,j}\,\omega_\nu\big({\mathfrak p}^{-j}x-u(k)\big),\eqno (5.4)$$

\parindent=0mm \vspace{.1in}
where $\{c_k^j\}_{k \in \mathbb N_0},\,\{d_k^{\nu,j}\}_{k \in \mathbb N_0} \in l^2(\mathbb N_0).$

\parindent=0mm \vspace{.2in}
Implementation of equation (5.1) for $n=0$ gives the decomposition of $f_{j}(x)$ as

$$\begin{array}{rcl}
f_j(x) &=& \displaystyle\sum_{k \in \mathbb N_0}c_k^j\omega_0\big({\mathfrak p}^{-j}x-u(k)\big)\\\\
&=&\dfrac{1}{\sqrt{q}}\displaystyle\sum_{k \in \mathbb N_0}c_k^j\displaystyle\sum_{\nu=0}^{q-1}\displaystyle\sum_{\mu \in \mathbb N_0}\tilde a_{k-q\mu}^\nu \omega_{qn+\nu}\big({\mathfrak p}^{-j}x-u(\mu)\big)\\\\
&=&\dfrac{1}{\sqrt{q}}\displaystyle\sum_{k \in \mathbb N_0}\displaystyle\sum_{\mu \in \mathbb N_0}c_\mu^j\displaystyle\sum_{\nu=0}^{q-1}\tilde a_{\mu-qk}^\nu \omega_\nu\big({\mathfrak p}^{-j+1}x-u(k)\big)\\\\
&=&\dfrac{1}{\sqrt{q}}\displaystyle\sum_{k \in \mathbb N_0}\left(\displaystyle\sum_{\mu \in \mathbb N_0}c_\mu^j\tilde a_{\mu-qk}^0\right) \omega_0\big({\mathfrak p}^{-j+1}x-u(k)\big)\\\\
&&\quad \quad \quad \quad \quad+ \;\dfrac{1}{\sqrt{q}}\displaystyle\sum_{k \in \mathbb N_0}\displaystyle\sum_{\nu=1}^{q-1}\left(\displaystyle\sum_{\mu \in \mathbb N_0}c_\mu^j\tilde a_{\mu-qk}^\nu \right)\omega_\nu\big({\mathfrak p}^{-j+1}x-u(k)\big)\\\\
&=&\displaystyle\sum_{k \in \mathbb N_0}c_k^{j-1}\varphi\big({\mathfrak p}^{-j+1}x-u(k)\big)+\displaystyle\sum_{\nu=1}^{q-1}\displaystyle\sum_{k \in \mathbb N_0}d_k^{i,j-1}\omega_\nu\big({\mathfrak p}^{-j+1}x-u(k)\big)\\\\
&=& f_{j-1}(x)+g_{j-1}(x),
\end{array}$$

\parindent=0mm \vspace{.1in}
where
$$  c_k^{j-1}=\dfrac{1}{\sqrt{q}}\displaystyle\sum_{\mu \in \mathbb N_0}c_\mu^j\,\tilde a_{\mu-qk}^0,\quad  d_k^{i,j-1}=\dfrac{1}{\sqrt{q}}\displaystyle\sum_{\mu \in \mathbb N_0}c_\mu^j\,\tilde a_{\mu-qk}^\nu,\eqno(5.5)$$

\parindent=0mm \vspace{.1in}
$k \in \mathbb N_0, j= J, J-1,\dots, J-M+1$.

\parindent=8mm \vspace{.2in}
For all $r \in \mathbb N_0,$ we have
$$g_j\in W_j=\Delta^jW_1=\Delta^{j-r}\Delta^rW_1=\Delta^{j-r}\bigoplus_{{\nu=q^r} }^{ q^{r+1}-1}E_\nu .$$

\parindent=0mm \vspace{.1in}
Using Theorem 5.1 for $n=1,2,\dots,q^{r+1}-1$,  yields

$$\begin{array}{rcl}
g_j(x)&=&\displaystyle\sum_{\nu=1}^{q-1}\displaystyle\sum_{k \in \mathbb N_0}d_k^{i,j}\omega_\nu\big({\mathfrak p}^{-j}x-u(k)\big)\\\\
&=&\dfrac{1}{\sqrt{q}}\displaystyle\sum_{k \in \mathbb N_0}\displaystyle\sum_{\nu=1}^{q-1}d_k^{i,j}\displaystyle\sum_{\mu \in \mathbb N_0}\displaystyle\sum_{s=o}^{q-1}\tilde a_{k-q\mu}^s \omega_{q\nu+\mu}\big({\mathfrak p}^{-j+1}x-u(\mu)\big)\\\\
&=&\dfrac{1}{\sqrt{q}}\displaystyle\sum_{k \in \mathbb N_0}\displaystyle\sum_{\nu=1}^{q^2-1}\left(\displaystyle\sum_{\mu \in \mathbb N_0}d_\mu^{\lfloor \nu/q \rfloor,j}\,\tilde a_{\mu-qk}^{\big(\nu-q\lfloor \nu/q \rfloor\big)} \right)\omega_\nu\big({\mathfrak p}^{-j+1}x-u(k)\big)\\\\
&=&\displaystyle\sum_{k \in \mathbb N_0}\displaystyle\sum_{\nu=1}^{q^2-1}d_k^{\nu, j, 1}\omega_\nu\big({\mathfrak p}^{-j+1}x-u(k)\big)\\
& \vdots&\\
&=&\displaystyle\sum_{k \in \mathbb N_0}\displaystyle\sum_{\nu=q^r}^{q^{r+1}-1}d_k^{\nu, j, r}\omega_\nu\big({\mathfrak p}^{-j+r}x-u(k)\big),
\end{array}$$

\parindent=0mm \vspace{.1in}
where
$$d_k^{\nu, j, i}=\dfrac{1}{\sqrt{q}}\displaystyle\sum_{\mu \in \mathbb N_0}d_\mu^{\lfloor \nu/q \rfloor,j,i-1}\,\tilde a_{\mu-qk}^{\big(\nu-q\lfloor \nu/q \rfloor\big)},\quad d_k^{\nu, j, 0}=d_k^{\nu, j}.\eqno (5.6)$$

\parindent=0mm \vspace{.1in}
$ i=1,2,\dots,r,\; \nu=q^i, q^i+1, \dots, q^{i+1}-1$.

\parindent=0mm \vspace{.2in}
Therefore, for $r \in \mathbb N_0, \, f_J$ can be decomposed as:
$$\begin{array}{rcl}
f_J&=&f_{J-M}+\displaystyle\sum_{j=J-M}^{J-1}g_j\\\
&=&\displaystyle\sum_{k \in \mathbb N_0}c_k^{J-M} \omega_0 \big({\mathfrak p}^{J-M}x-u(k)\big)+\displaystyle\sum_{j=J-M}^{J-1}\displaystyle\sum_{\nu=1}^{q-1}\sum_{k \in \mathbb N_0}d_k^{\nu,j}\omega_\nu\big({\mathfrak p}^{-j}x-u(k)\big)\\\
&=&\displaystyle\sum_{k \in \mathbb N_0}c_k^{J-M} \varphi \big({\mathfrak p}^{J-M}x-u(k)\big)+\displaystyle\sum_{j=J-M}^{J-1}\displaystyle\sum_{\nu=q^r}^{q^{r+1}-1}\sum_{k \in \mathbb N_0}d_k^{\nu,j,r}\omega_\nu\big({\mathfrak p}^{r-j}x-u(k)\big),
\end{array} $$

\parindent=0mm \vspace{.1in}
where the coefficients are given by the equations (5.5) and (5.6).

\parindent=0mm \vspace{.2in}
On the other hand, by using equation (3.1), we can reconstruct $g_j(\cdot)$ as follows:
$$\begin{array}{rcl}
g_j(x)&=&\displaystyle\sum_{\nu=q^r}^{q^{r+1}-1}\sum_{k \in \mathbb N_0}d_k^{\nu,j,r}\omega_\nu\big({\mathfrak p}^{r-j}x-u(k)\big)\\\\
&=&\displaystyle\sum_{\nu=q^r}^{q^{r+1}-1}\sum_{k \in \mathbb N_0}d_k^{\nu,j,r}\displaystyle\sum_{\mu \in \mathbb N_0}a_\mu^{\big(\nu-q\lfloor \nu/q \rfloor\big)}\omega_{\lfloor \nu/q \rfloor}\big({\mathfrak p}^{r-j-1}x-u(qk-\mu)\big)\\\\
&=&\displaystyle\sum_{\nu=q^{r-1}}^{q^r-1}\sum_{k \in \mathbb N_0}d_k^{\nu,j,r-1}\omega_\nu({\mathfrak p}^{r-j-1}x-u(k)\big)\\\\
&\vdots&\\\\
&=&\displaystyle\sum_{\nu=1}^{q-1}\sum_{k \in \mathbb N_0}d_k^{\nu,j}\omega_\nu({\mathfrak p}^{-j}x-u(k)\big),
\end{array} $$

\parindent=0mm \vspace{.1in}
where
$$d_k^{\nu, j, i-1}=\displaystyle\sum_{s=0}^{q-1}\displaystyle\sum_{\mu \in \mathbb N_0}d_\mu^{q\nu+s,j,i}\, a_{k-q\mu}^s,\quad d_k^{\nu, j}=d_k^{\nu, j, 0}.
 \eqno (5.7)$$

\parindent=0mm \vspace{.1in}
 $i=1,2,\dots,r,~  \nu=q^{i-1}, q^{i-1}+1, \dots, q^{i}-1.$

\parindent=0mm \vspace{.2in}
Thus, after obtaining the coefficients $d_k^{\nu,j},\,\nu=1,2,\dots,q-1,\, j=J-M,\dots,J-1,\,k \in \mathbb N_0$, we use Theorem 5.1 and (2.6) to construct $f_j$ as follows:

$$\begin{array}{rcl}
f_j&=&f_{j-1}+g_{j-1}\\\\
&=&\displaystyle\sum_{k \in \mathbb N_0}c_k^{j-1} \omega_0 \big({\mathfrak p}^{-j-1}x-u(k)\big)+\displaystyle\sum_{\nu=1}^{q-1}\sum_{k \in \mathbb N_0}d_k^{\nu,j-1}\omega_\nu\big({\mathfrak p}^{-j+1}x-u(k)\big)\\\\
&=&\displaystyle\sum_{k \in \mathbb N_0}c_k^{j-1}\displaystyle\sum_{\mu \in \mathbb N_0}a_\mu^0 \omega_0 \big({\mathfrak p}^{-j}x-u(qk-\mu)\big)+\displaystyle\sum_{\nu=1}^{q-1}\sum_{k \in \mathbb N_0}d_k^{\nu,j-1}\displaystyle\sum_{\mu \in \mathbb N_0}a_\mu^\nu \omega_0\big({\mathfrak p}^{-j}x-u(qk-\mu)\big)\\\\
&=&\displaystyle\sum_{k \in \mathbb N_0}c_k^{j-1}\displaystyle\sum_{\mu \in \mathbb N_0}a_{\mu-qk}^0 \omega_0 \big({\mathfrak p}^{-j}x-u(\mu)\big)+\displaystyle\sum_{\nu=1}^{q-1}\sum_{k \in \mathbb N_0}d_k^{\nu,j-1}\displaystyle\sum_{\mu \in \mathbb N_0}a_{\mu-qk}^\nu \omega_0\big({\mathfrak p}^{-j}x-u(\mu)\big)\\\\
&=&\displaystyle\sum_{k \in \mathbb N_0}\left(\displaystyle\sum_{\mu \in \mathbb N_0}c_k^{j-1}a_{\mu-qk}^0 +\displaystyle\sum_{\nu=1}^{q-1}\sum_{\mu \in \mathbb N_0}d_\mu^{\nu,j-1}a_{\mu-qk}^\nu\right) \omega_0\big({\mathfrak p}^{-j}x-u(k)\big)
\end{array} $$

\parindent=0mm \vspace{.2in}
\qquad $=\displaystyle\sum_{k \in \mathbb N_0}c_k^j\varphi\big({\mathfrak p}^{-j}x-u(k)\big),$

\parindent=0mm \vspace{.1in}
where
$$c_k^j=\displaystyle\sum_{\mu \in \mathbb N_0}c_k^{j-1}a_{\mu-qk}^0 +\displaystyle\sum_{\nu=1}^{q-1}\sum_{\mu \in \mathbb N_0}d_\mu^{\nu,j-1}a_{\mu-qk}^\nu,\; j=J-M+1, J-M+2,\dots,J,\;k \in \mathbb N_0.\eqno(5.8)$$

Therefore, with the given sequences $\left\{c_k^{J-M}\right\}_{k \in \mathbb N_0}$ and $\left\{d_k^{\nu,J-M}\right\}_{k \in \mathbb N_0}, \nu=1,2,\dots,q-1 $, and applying (5.8), one can reconstruct

$$f\approx f_J= \sum_{k \in \mathbb N_0}c_k^J \omega_0\big({\mathfrak p}^{-J}x-u(k)\big)\in V_J.$$

\parindent=0mm \vspace{.2in}

{\bf{References}}

\begin{enumerate}

{\small{
\item  B. Behera and Q. Jahan, Biorthogonal wavelets on local fields of positive characteristic, Comm. Math. Anal. 11 (2013), 52-75.

\item  B. Behera and Q. Jahan, Wavelet packets and wavelet frame packets on local fields of positive characteristic,, J. Math. Anal. Appl. 395 (2012), 1-14.

\item J. J. Benedetto and R. L. Benedetto,  A wavelet theory for local fields and related groups, J. Geom. Anal. 14 (2004), 423-456.

\item Q. Chen and Z. Chang,  A study on compactly supported orthogonal vector-valued wavelets and wavelet packets,  Chaos, Solitons and Fractals, {31} (2007), 1024-1034.

\item C. Chui  and  C. Li, Non-orthogonal wavelet packets, SIAM J. Math. Anal.  24(3) (1993), 712-738.

\item A. Cohen and I. Daubechies, On the instability of arbitrary biorthogonal wavelet packets,  SIAM J. Math. Anal. 24(5) (1993), 1340-1354.

\item R. R. Coifman, Y. Meyer, S. Quake and  M. V. Wickerhauser, Signal processing and compression with wavelet packets,  Technical Report, Yale University, 1990.

\item Yu. A. Farkov,  Orthogonal wavelets with compact support on locally compact Abelian groups,   Izv. Math. 69(3) (2005), 623-650.

\item H. K. Jiang,  D. F. Li and  N. Jin, Multiresolution analysis on local fields,  J. Math. Anal. Appl. 294 (2004), 523-532.

\item A. Yu. Khrennikov, V. M. Shelkovich and M. Skopina, $p$-Adic refinable functions and MRA-based wavelets, J. Approx. Theory. 161 (2009),  226-238.

\item W. C. Lang, Orthogonal wavelets on the Cantor dyadic group,  SIAM J. Math. Anal. 27 (1996), 305-312.

\item J. N. Leng, Z. Cheng, T. Huang and C. Lai, Construction and properties of multiwavelet packets with arbitrary scale and the related algorithms of decomposition and reconstruction, Comput. Math. Appl. 51 (2006), 1663-1676.

\item D. F. Li and H. K. Jiang, The necessary condition and sufficient conditions for wavelet frame on local fields,  J. Math. Anal. Appl. 345(2008), 500-510.

\item S. F. Lukomskii, Multiresolution analysis on product of zero-dimensional Abelian groups, J. Math. Anal. Appl. 385 (2012), 1162-1178.

\item D. Ramakrishnan and  R. J. Valenza,  {\it Fourier Analysis on Number Fields}, Graduate Texts in Mathematics 186, Springer-Verlag, New York, 1999.

\item F. A. Shah,  Construction of wavelet packets on $p$-adic field,  Int. J. Wavelets Multiresolut. Inf. Process. 7(5) (2009), 553-565.

\item F. A. Shah, Biorthogonal $p$-wavelet packets related to the Walsh polynomials,  J. Classical Anal. 2 (2013).

\item F. A. Shah and L. Debnath,  Tight wavelet frames on local fields, Analysis 33(2013), 293-307.

\item Z. Shen, Non-tensor product wavelet packets in $L^2(\mathbb R^s)$, SIAM J. Math. Anal. 26(4) (1995), 1061-1074.

\item T. Stavropoulos  and M. Papadakis,  On the multiresolution analyses of abstract Hilbert spaces,  Bull. Greek Math. Soc.  40 (1998), 79-92.

\item M. H. Taibleson,  {\it{Fourier Analysis on Local Fields,}} Princeton University Press, 1975.
}}

\end{enumerate}

\end{document}